\documentclass[12pt]{amsart}
\usepackage{amsmath, amssymb,graphicx,caption, subcaption, mathtools, enumitem}
\usepackage[pdftex,hyperindex=true]{hyperref}
\usepackage{url}
\usepackage{etoolbox}
\patchcmd{\section}{\scshape}{\bfseries}{}{}
\patchcmd{\subsection}{\bfseries}{\itshape}{}{}

\makeatletter
\def\@seccntformat#1{%
  \protect\textup{\protect\@secnumfont
    \ifnum\pdfstrcmp{section}{#1}=0 \bfseries\fi
    \ifnum\pdfstrcmp{subsection}{#1}=0 \itshape \fi
    \csname the#1\endcsname
    \protect\@secnumpunct
  }%
}  
\makeatother

\usepackage[utf8]{inputenc}
\usepackage[english]{babel}

\DeclarePairedDelimiter\abs{\lvert}{\rvert}%
\DeclarePairedDelimiter\norm{\lVert}{\rVert}%
\theoremstyle{plain}

\newtheorem{conjecture}{Conjecture}[section]

\theoremstyle{definition}

\newtheoremstyle{note}
{3pt}
{3pt}
{\itshape}
{}
{\itshape}
{:}
{.5em}
{}

\theoremstyle{note}

\theoremstyle{plain} 
\newcommand{\thistheoremname}{}
\newtheorem*{genericthm}{\thistheoremname}

\numberwithin{equation}{section}

\newcommand{\acknowledge}{\subsection*{Acknowledgements}}

\def\Dbar{\leavevmode\lower.6ex\hbox to 0pt{\hskip-.23ex \accent"16\hss}D}

\begin{document}

\title{The Atiyah-Sutcliffe determinant}
\author{J. Malkoun}
\address{Department of Mathematics and Statistics\\
Faculty of Natural and Applied Sciences\\
Notre Dame University-Louaize, Zouk Mosbeh\\
P.O.Box: 72, Zouk Mikael,
Lebanon}
\email{joseph.malkoun@ndu.edu.lb}

\date{March 14, 2019}

\maketitle

\begin{abstract} We present a general formula for the Atiyah-Sutcliffe determinant function, which holds for any integer $n \geq 2$, as 
a global factor times a sum of terms, with each term similar to a higher degree cross-ratio. The formula is to our knowledge new. 

We also conjecture that the Atiyah-Sutcliffe determinant is a rational 
linear combination of products of factors of only two simple types, each of them manifestly $SO(3)$-invariant. This allows us to obtain a 
conjectural purely angular formula for the determinant for $n=4$, as an illustration of how our conjecture can be applied.
\end{abstract}

\section{Introduction} \label{intro}

Denoting by $C_n(\mathbb{R}^3)$ the configuration space of $n$ distinct points in $\mathbb{R}^3$, we consider a configuration
\[ \mathbf{x} = (\mathbf{x}_1,\ldots,\mathbf{x}_n) \in C_n( \mathbb{R}^3 ) .\]
Given $a \neq b$, with $1 \leq a,b \leq n$, we form
\[ v_{ab} = \frac{\mathbf{x}_b - \mathbf{x}_a}{\norm{\mathbf{x}_b - \mathbf{x}_a}} \in S^2 ,\] 
where $\norm{.}$ denotes the Euclidean norm in $\mathbb{R}^3$. Stereographic projection is a smooth map, which we will denote by $s$, from $S^2$ onto the Riemann 
sphere $\hat{\mathbb{C}}$, and can be defined as follows:
\[ s(x,y,z) = \begin{cases} \frac{x+iy}{1-z} \text{ , if $(x,y,z) \neq (0,0,1)$}\\
\infty  \text{ , otherwise.} \end{cases}\]
We then define, for each $a \neq b$ ($1 \leq a,b \leq n$),
\[ t_{ab} = s(v_{ab}) .\]
We now form, for each $a$, $1 \leq a \leq n$, the polynomial
\[ p_a(t) = \prod_{b \neq a} (t - t_{ab}) ,\]
with the understanding that factors $t-t_{ab}$ corresponding to $t_{ab} = \infty$  are set to $1$ (this can be justified using homogeneous coordinates on $P^1(\mathbb{C})$). 

The Atiyah-Sutcliffe conjecture $1$ can now be formulated.

\begin{conjecture}[AS conjecture $1$]Given any $\mathbf{x} \in C_n(\mathbb{R}^3)$, the corresponding $n$ polynomials $p_a$, for $1 \leq a \leq n$, are linearly independent over $\mathbb{C}$.
\end{conjecture}

We now define the Atiyah-Sutcliffe normalized determinant function $D: C_n(\mathbb{R}^3) \to \mathbb{C}$.
Let
\[ \delta(\mathbf{t}) = \prod_{1 \leq a<b \leq n} (t_{ab} - t_{ba}) ,\]
where $\mathbf{t} = (t_{ab})$. I am using inhomogeneous coordinates $t_{ab}$ on the affine subset $\mathbb{C}$ of $P^1(\mathbb{C})$ corresponding to $u \neq 0$, where $[u:v]$ is a 
homogeneous coordinate of a point on $P^1(\mathbb{C})$. However, it is straightforward to figure out what each factor means if some of the $t_{ab}$ are infinity; indeed, 
it suffices to use homogeneous coordinates and $2$-by-$2$ determinants, instead of inhomogeneous coordinates and differences.
Form the matrix
\[ A = (p_1,\ldots,p_n) \]
having the coefficients of $p_j$ corresponding to decreasing powers of $t$ as its $j$-th column.

We can now define the normalized Atiyah-Sutcliffe determinant $D$ by
\[ D(\mathbf{x}) = \frac{\det(A)}{\delta(\mathbf{t})} \]

The Atiyah-Sutcliffe conjecture $2$ can now be stated.

\begin{conjecture}[AS conjecture $2$]Given any $\mathbf{x} \in C_n(\mathbb{R}^3)$, 
\[\abs{D(\mathbf{x})} \geq 1 .\]
\end{conjecture}

In section \ref{general}, we formulate and prove a general formula for the Atiyah-Sutcliffe normalized determinant $D$. In section \ref{conjecture}, 
we conjecture that $D$ can be expressed as a rational linear combination of terms containing factors of only two types:
\[ (v_{ab}, v_{cd}) \text{ and } \det(v_{ab}, v_{cd}, v_{ef}) ,\]
and then illustrate the usefulness of this conjecture by providing a conjectural purely angular formula for $D$ for $n=4$ (Conjecture \ref{angular}).

\section{A general formula for the Atiyah-Sutcliffe normalized determinant} \label{general}

Define
\[ I_n = \{ (a,b) ; 1 \leq a,b \leq n \text{ and } a \neq b \} \]
and let $G$ be the group of all permutations of $I_n$ which fix the first coordinate, namely, for every $\sigma \in G$ and every $(a,b) \in I_n$, we have
\[ \sigma (a,b) = (a, b') \]
for some $b' \neq a$ ($1 \leq b' \leq n$). We note that the order of $G$ is
\[ \abs{G} =  \left( (n-1)! \right)^n .\]
The action of $G$ on $I_n$ induces an action of $G$ on the polynomial space $\mathbb{C}[\mathbf{t}]$. For instance, we have
\[ (\sigma.\delta)(\mathbf{t}) = \prod_{1 \leq a<b \leq n} (t_{\sigma(a,b)} - t_{\sigma(b,a)}) \]
Given an integer $n \geq 2$, we introduce the integer
\[ c_n = \prod_{a = 1}^{n-1} (a!)^2 \]
Our formula for the Atiyah-Sutcliffe determinant can now be formulated.
\[ D(\mathbf{x}) = \frac{1}{c_n} \sum_{\sigma \in G} \frac{(\sigma.\delta)(\mathbf{t})}{\delta(\mathbf{t})} \]

The idea of the proof is that the subspace of $\mathbb{C}[\mathbf{t}]$ having the same symmetries as the numerator of the Atiyah-Sutcliffe determinant is 
$1$-dimensional. More precisely, let $W$ be the subspace of $\mathbb{C}[\mathbf{t}]$ consisting of polynomials which are symmetric under the action of $G$, multilinear and 
skew-symmetric as a function of the $t_{ab}$, thought of as elements of $\mathbb{C}^2 \setminus \{ \mathbf{0} \}$ 
up to scaling by a complex factor (i.e. using homogeneous coordinates on $\mathbf{P}^1(\mathbb{C})$). Then $W$ is 
complex $1$-dimensional.

Moreover, one can prove that
\[ \sum_{\sigma \in G} (\sigma.\delta)(\mathbf{t}) \]
(essentially the numerator on the right-hand side of our formula) is skew-symmetric under the action of the symmetric group $\Sigma_n$. It remains to calculate the 
global factor, which can be done by considering for instance collinear configurations, at which it is known that $D$ takes the value $1$. This proves our general formula.

\section{A general conjecture for $D$ and an application} \label{conjecture}

We make the following conjecture.

\begin{conjecture}The Atiyah-Sutcliffe normalized determinant $D$ can be expressed as a rational linear combination of terms containing factors of only two types:
\[ (v_{ab}, v_{cd}) \text{ and } \det(v_{ab}, v_{cd}, v_{ef}) ,\]
where the $v_{ab}$ are as defined in section \ref{intro}. \label{main_conj}
\end{conjecture}

We illustrate the usefulness of this conjecture by providing a conjectural purely angular formula for $D$ for $n=4$.

\begin{conjecture}If $n=4$, then
\begin{align*} \operatorname{Re}(D) &= \frac{3}{8} + \frac{3}{2} \operatorname{Av}((v_{12},v_{13})) +  \frac{3}{2} \operatorname{Av}((v_{12},v_{13}) (v_{14},v_{24})) + \cdots \\ 
   &\cdots + \frac{3}{8} \operatorname{Av}((v_{12},v_{34}) (v_{13},v_{24})) + \frac{1}{2} \operatorname{Av}((v_{12},v_{14}) (v_{13},v_{23}) (v_{24},v_{34})) \end{align*}
and
\begin{align*} \operatorname{Im}(D) &= -\frac{1}{32} \operatorname{Av}(\det(v_{12},v_{13},v_{24})) - \frac{1}{32} \operatorname{Av}(\det(v_{12},v_{14},v_{23}) (v_{24},v_{34}))
\end{align*}
where $\operatorname{Av}$ denotes averaging over the symmetric group $\Sigma_4$. \label{angular}
\end{conjecture}

The way we used conjecture \ref{main_conj} in order to obtain the conjectural angular formula for $n=4$ (i.e. conjecture \ref{angular}) is as follows. We used a computer to 
generate all possible terms of the right degree which are of the form given in conjecture \ref{main_conj}, and formed the distinct $\Sigma_4$ orbits of such terms. Then using randomly 
generated configurations, whose number is equal to the number of unknown coefficients, a linear system was set up for these coefficients, and was then very easily solved by the computer.

The author hopes that such angular formulas will turn out to be useful. They are to his knowledge new. The reader may wish to also consult for instance the article \cite{EN} by Eastwood and Norbury, 
who actually prove a formula for $D$ for $n=4$, but their formula is not (explicitly) angular.

\acknowledge{The author is forever indebted to Sir Michael Atiyah for many long discussions on this problem. Sir Michael did not actually believe that using formulas was the right way to approach this problem. 
I remember very well that when I once tried to explain to him my formulas, he interrupted me and told me that a solution should be geometric and beautiful! However, I choose to make my formulas available, 
hoping someone may make good use of them. The mathematical genius Sir Michael will be greatly missed, but he did leave behind a wealth of wisdom, and a large number of beautiful results.}

\pdfbookmark[1]{References}{ref}

\end{document}